\documentclass{amsart}
\usepackage{amsmath,amsfonts,amssymb,amsthm}
\DeclareMathOperator{\GL}{\text{GL}}

\DeclareMathOperator{\lcm}{\text{lcm}}
\DeclareMathOperator{\Gal}{\text{Gal}}
\newcommand{\BR}{{\mathbb R}}
\newcommand{\BQ}{{\mathbb Q}}
\newcommand{\BC}{{\mathbb C}}
\newcommand{\CO}{{\mathcal O}}
\newcommand{\rc}{{\check r}}
\newcommand{\inv}{{^{-1}}}
\newtheorem*{theorem}{Theorem}
\newtheorem*{definition}{Definition}
\newcommand{\lexp}[2]{\kern\scriptspace\vphantom{#2}^{#1}\kern-\scriptspace#2}
\begin{document}
\author{J.~Michel}
\title[Hurwitz action]{Hurwitz action on tuples of Euclidean reflections}
\date{10th August, 2004}
\maketitle
This  note  was prompted  by  the  reading of  \cite{Humphries},  which
purports to  show that if an  $n$-tuple of Euclidean reflections  has a
finite orbit  under the  Hurwitz action  of the  braid group,  then the
generated group  is finite. I noticed  that the proof given  is fatally
flawed \footnote{The problem is in  proposition 2.3, which is essential
to the  main theorem (1.1)  of the paper.  The argument given  there is
basically that if a Coxeter group has a reflection representation where
the image of the Coxeter element is  of finite order, then the image of
that representation is finite. However this is false: the Cartan matrix
$\bigl(\begin{smallmatrix}
2&-1&-1\\-1&2&-l\\-1&-l&2\\
\end{smallmatrix}\bigr)$ where $l^2+l=\sqrt 2$
defines  Euclidean  reflections  which  give  a  representation  of  an
infinite rank 3 Coxeter group, such that the image of the Coxeter group
is  infinite  but the  image  of  the Coxeter  element  is  of order  8
(personal communication  of F.Zara).};  however, using the  argument of
Vinberg given in  \cite{DM}, I found a short  (hopefully correct) proof
which  at  the  same  time considerably  simplifies  the  computational
argument given  in \cite{DM}.  This is  what I  expound below.  I first
recall  all the  necessary  notation and  assumptions, expounding  some
facts in slightly more generality than necessary.

\subsection{Hurwitz action}
\begin{definition}
Given a group $G$, we call  Hurwitz action the action of the $n$-strand
braid group $B_n$ with standard generators $\sigma_i$ on $G^n$ given by
$$\sigma_i(s_1,\ldots,s_n)=(s_1,\ldots,s_{i-1},s_{i+1},s_i^{s_{i+1}},
s_{i+2},\ldots,s_n).$$
\end{definition}
The  inverse  is given  by  $\sigma_i^{-1}(s_1,\ldots,s_n)=(s_1,\ldots,
s_{i-1},  \lexp{s_i}s_{i+1},s_i,s_{i+2},\ldots,s_n)$.   Here  $a^b$  is
$b^\inv ab$ and $\lexp ba$ is $bab\inv$.

This  action  preserves  the  product  of  the  $n$-tuple.  We  need  to
repeat  some  remarks  in   \cite{Humphries}.  By  decreasing  induction
on   $i$   one    sees   that   $\sigma_i\ldots\sigma_n(s_1,\ldots,s_n)=
(s_1,\ldots,s_{i-1},s_n,s_i^{s_n},\ldots,s_{n-1}^{s_n})$.  In particular
if  $\gamma=\sigma_1\ldots\sigma_{n-1}$  we get  $\gamma(s_1,\ldots,s_n)
=(s_n,s_1,\ldots,s_{n-1})^{s_n}$  whence, if  $c=s_1\ldots s_n$,  we get
that $\gamma^n(s_1,\ldots,s_n)=(s_1,\ldots,s_n)^c$.

We  also  deduce  that  given  any  subsequence  $(i_1,\ldots,i_k)$  of
$(1,\ldots,n)$,  there  exists  an  element of  the  Hurwitz  orbit  of
$(s_1,\dots,s_n)$ which begins by $(s_{i_1},\ldots,s_{i_k})$.

Assume now that the Hurwitz orbit of $(s_1,\ldots,s_n)$ is finite. Then
some power of $\gamma$ fixes $(s_1,\ldots,s_n)$, thus some power of $c$
is  central in  the  subgroup  generated by  the  $s_i$. Similarly,  by
looking at the action of  $\sigma_1\ldots\sigma_{k-1}$ on an element of
the orbit beginning  by $(s_{i_1},\ldots,s_{i_k})$ we get  that for any
subsequence $(i_1,\ldots,i_k)$  of $(1,\ldots,n)$ there exists  a power
of  $s_{i_1}\ldots  s_{i_k}$  central  in  the  subgroup  generated  by
$(s_{i_1},\ldots,s_{i_k})$.

\subsection{Reflections}

Let $V$ be a  vector space on some subfield $K$ of  $\BC$. We call {\it 
complex  reflection} a  finite order  element $s\in\GL(V)$  whose fixed 
points are  a hyperplane. If  $\zeta$ (a root  of unity) is  the unique 
non-trivial  eigenvalue  of $s$,  the  action  of  $s$ can  be  written 
$s(x)=x-\rc(x)r$ where $r\in V$ and $\rc$  is an element of the dual of 
$V$  satisfying  $\rc(r)=1-\zeta$.  These  elements are  unique  up  to 
multiplying $r$  by a scalar  and $\rc$ by  the inverse scalar.  We say 
that $r$ (resp. $\rc$) is a root (resp. coroot) associated to $s$.      

\subsection{Cartan Matrix}
If  $(s_1,\ldots,s_n)$  is  a  tuple  of  complex  reflections  and  if
$r_i,\rc_i$ are  corresponding roots and  coroots, we call  {\it Cartan
matrix} the  matrix $C=\{\rc_i(r_j)\}_{i,j}$. This matrix  is unique up
to conjugating  by a  diagonal matrix. Conversely,  a class  modulo the
action of diagonal  matrices of Cartan matrices is an  invariant of the
$\GL(V)$-conjugacy class of  the tuple. It determines this  class if it
is invertible and $n=\dim V$. Indeed,  this implies that the $r_i$ form
a basis  of $V$; and  in this basis the  matrix $s_i$ differs  from the
identity  matrix only  on the  $i$-th line,  where the  opposed of  the
$i$-th line of $C$ has been added; thus $C$ determines the $s_i$.

 If $C$  can be chosen  Hermitian (resp.  symmetric), such a  choice is
then unique up to conjugating by a diagonal matrix of norm $1$ elements
of $K$ (resp. of signs).

 If $C$  is Hermitian (which  implies that the  $s_i$ are of  order 2), 
then the sesquilinear  form given by ${}^tC$ is invariant  by the $s_i$ 
(if  the  $s_i$ are  not  of  order $2$,  but  the  matrix obtained  by 
replacing  all  elements  on  the  diagonal  of  ${}^tC$  by  $2$'s  is 
Hermitian, then the latter matrix defines a sesquilinear form invariant 
by the $s_i$).                                                          

\subsection{Coxeter element}
We keep the notation as above and we assume that the $r_i$ form a basis 
of $V$. We recall a result of \cite{Coleman} on the ``Coxeter'' element 
$c=s_1\ldots s_k$.  If we write  $C=U+V$ where $U$ is  upper triangular 
unipotent  and  where $V$  is  lower  triangular (with  diagonal  terms 
$-\zeta_i$, thus $V$ is also unipotent when $s_i$ are of order 2), then 
the matrix of $c$  in the $r_i$ basis is $-U\inv V$  (to see this write 
it as $Us_1\ldots s_n=-V$ and look at partial products in the left-hand 
side starting  from the left).  As $U$ is  of determinant 1,  we deduce 
that $\chi(c)=\det(xI+U\inv V)=\det(xU+V)$  where $\chi(c)$ denotes the 
characteristic polynomial ;  in particular $\det(C)=\chi(c)\mid_{x=1}$; 
one also gets that the fix-point set of $c$ is the kernel of $C$, equal 
to the intersection of the reflecting hyperplanes.                      

\subsection{The main theorem}
The  next  theorem  implies  the statement  given  in  \cite{Humphries} 
(\cite[1.1]{Humphries} considers  Euclidean reflections with  the $r_i$ 
linearly independent; if  the $r_i$ are chosen of the  same length this 
implies that $C$  is symmetric, and as  $C$ is then the  Gram matrix of 
the $r_i$ it is invertible):                                            
\begin{theorem}
Let $(s_1,\ldots,s_n)$ be  a tuple of reflections  in $\GL(\BR^n)$ which
have  an associated  Cartan matrix  symmetric and invertible. Assume  in
addition that the  Hurwitz orbit of the tuple is  finite. Then the group
generated by  the $s_i$ is finite.
\end{theorem}
\begin{proof}   In   the   next    paragraph,   we   just   need   that 
$(s_1,\ldots,s_n)$  is a  tuple of  complex reflections  with a  finite 
Hurwitz orbit and with the $r_i$ a basis of $V$.                        

A  straightforward  computation  shows  than  an  element  of  $\GL(V)$ 
commutes  to  the  $s_i$  if  and  only if  it  acts  as  a  scalar  on 
the  subspaces  generated  by  $\{  r_i\}_{i\in  I}$  where  $I$  is  a 
block  of  $C$   (i.e.,  a  connected  component  of   the  graph  with 
vertices $\{1,\ldots,n\}$  and edges  $(i,j)$ for  each pair  such that 
either  $C_{i,j}$ or  $C_{j,i}$ is  not  zero). The  finiteness of  the 
Hurwitz orbit  implies that  for any subsequence  $(i_1,\ldots,i_k)$ of 
$(1,\ldots,n)$, there  exists a power of  $s_{i_1}\ldots s_{i_k}$ which 
commutes to $s_{i_1},\ldots,s_{i_k}$. This power  acts thus as a scalar 
on each subspace generated by the $r_{i_j}$ in a block of the submatrix 
of $C$  determined by  $(i_1,\ldots,i_k)$. As  the determinant  of each 
$s_{i_j}$ on  this subspace is  a root of unity,  the scalar must  be a 
root of unity. Thus, the restriction of each $s_{i_1}\ldots s_{i_k}$ to 
the subspace  $<r_{i_1},\ldots,r_{i_k}>$ generated by the  $r_{i_j}$ is 
of finite order.                                                        

We use from now  on all the assumptions of the  theorem. Thus the $s_i$
are order $2$ elements of $O(C)$, the orthogonal group of the quadratic
form defined by $C$.

Also, $\chi(c)$  is a polynomial  with real  coefficients. As $c$  is of
finite order,  any real root of  $\chi(c)$ is $1$ or  $-1$. This implies
that $\chi(c)\mid_{x=1}$ is a nonnegative real number, and thus $\det C$
also.  The same  holds for any  principal  minor of  $C$,  since such  a
minor is  $\chi(c')\mid _{x=1}$  where $c'$ is  the restriction  of some
$s_{i_1}\ldots  s_{i_k}$ to  $<r_{i_1},\ldots, r_{i_k}>$.  The quadratic
form  defined by  $C$ is  thus  positive, and  as  $\det C\ne  0$ it  is
positive definite (cf. \cite[\S 7, exercice 2]{Bourbaki}).

We  now  digress  about  the  Cartan matrix  of  two  reflections  $s_1$
et  $s_2$.   Such  a  matrix   is  of  the   form  $\begin{pmatrix}2&a\\
b&2\\\end{pmatrix}$. If  $a=0$ and $b\ne 0$  or $a\ne 0$ and  $b=0$ then
$s_1s_2$  is  of  infinite  order.  Otherwise,  the  number  $ab$  is  a
complete invariant of  the conjugacy class of  $(s_1,s_2)$ restricted to
$<r_1,r_2>$, and $s_1s_2$ restricted to this subspace is of finite order
$m$ if and only  if there exists $k$ prime to  $m$ such that $ab=4\cos^2
k\pi/m$.

Since  $C$  is  symmetric  and  since the  restriction  of  $s_is_j$  to
$<r_i,r_j>$  is  of  finite  order, there  exists  prime  integer  pairs
$(k_{i,j},m_{i,j})$ such that $C_{i,j}=\pm 2 \cos k_{i,j}\pi/m_{i,j}$. If
$K$ is the cyclotomic subfield containing the $\lcm(2m_{i,j})$-th roots
of unity,  and if  $\CO$ is  the ring of  integers of  $K$, we  get that
all  coefficients  of $C$  lie  in  $\CO$. It  follows,  if  $G$ is  the
group  generated  by  the  $s_i$,  that  in  the  $r_i$  basis  we  have
$G\subset\GL(\CO^n)$.

We   now  apply   Vinberg's  argument   as  in   \cite[1.4.2]{DM}.  Let 
$\sigma\in\Gal(K/\BQ)$.  Then $\sigma(C)$  is again  positive definite: 
all arguments used  to prove that $C$ is positive  definite still apply 
for  $\sigma(C)$: it  is real,  symmetric, invertible  and the  Hurwitz 
orbit  of  $(\sigma(s_1),\ldots,\sigma(s_n))$  is still  finite.  Since 
$G\subset O(C)$, which  is compact, the entries of the  elements of $G$ 
in the  $r_i$ basis are of  bounded norm. Since $O(\sigma(C))$  is also 
compact for any $\sigma\in\Gal(K/\BQ)$, we get that entries of elements 
of $G$  are elements of  $\CO$ all of  whose complex conjugates  have a 
bounded norm.  There is  a finite  number of such  elements, so  $G$ is 
finite.                                                                 
\end{proof}


\begin{thebibliography}{99}
\bibitem[1]{Bourbaki}  N.Bourbaki, ``Alg\`ebre'', Chap. 9, {\sl Hermann
}, 1959.
\bibitem[2]{Coleman}  A.J.Coleman, ``Killing and the Coxeter transformation
of Kac-Moody algebras'', {\sl Invent. Math. \bf 95}(1989) 447--478.
\bibitem[3]{DM}  B.Dubrovin and  M.Mazzocco, ``Monodromy
of  certain Painlev\'e-VI  transcendents and  reflection groups'',  {\sl
Invent. Math. \bf 141}(2000), 55--147.
\bibitem[4]{Humphries}  S.~P.~Humphries,  ``Finite Hurwitz  braid  group
actions on  sequences of  Euclidean reflections'',  {\sl J.  Algebra \bf
269} (2003), 556--588.
\end{thebibliography}
\end{document}